\documentclass[12pt]{amsart}
\newtheorem{theorem}{Theorem}%[section]
\newtheorem{corollary}{Corollary}%[section]
\newtheorem{lemma}{Lemma}%[section]
\newtheorem{proposition}{Proposition}%[section]

\newtheorem{conjecture}{Conjecture}
\newtheorem{definition}{Definition}
\usepackage{amsmath,amstext,amssymb,amsopn,amsthm}
\usepackage{amsmath,amssymb,amsthm}
\usepackage[mathscr]{eucal}
\newtheorem{question}{Question}
\newtheorem{remark}{Remark}
\numberwithin{equation}{section}
\numberwithin{conjecture}{section}
\numberwithin{question}{section}
\numberwithin{definition}{section}
\numberwithin{theorem}{section}
\numberwithin{remark}{section}
\numberwithin{lemma}{section}
\numberwithin{proposition}{section}
\numberwithin{corollary}{section}

\begin{document}
\newcommand{\bD}{\mathrm{I\! D\!}}
\newcommand{\calD}{\mathcal{D}}
\newcommand{\calC}{\mathcal{C}}
\newcommand{\calE}{\mathcal{E}}
\newcommand{\calF}{\mathcal{F}}
\newcommand{\Rd}{\mathbb{R}^d}
\newcommand{\Rt}{\mathbb{R}^2}
\newcommand{\BR}{\mathcal{B}(\Rd)}
\newcommand{\R}{\mathbb{R}}
\newcommand{\C}{\mathbb{C}}
\newcommand{\al}{\alpha}
\newcommand{\ga}{\gamma}
\newcommand{\om}{\omega}
\newcommand{\G}{\Gamma}
\def\rayo{\leftarrow}
\def\tr{\triangle}
\newtheorem{cl}{Claim}
\newcommand{\bR}{\mathbf{R}}
\newtheorem{cor}{Corollary}
\theoremstyle{definition}
\newtheorem{fact}{Fact} \renewcommand{\thefact}{}
%\begin{document}

\title[Ground State Eigenfunction]{On the shape of the ground state eigenfunction for stable
processes}
\author{Rodrigo Ba\~{n}uelos}
\address{Mathematics Department, Purdue University, West Lafayette, IN 47907}
\email{banuelos@math.purdue.edu}
\thanks{R. Ba\~nuelos was supported in part by NSF grant
\# 9700585-DMS }
\author{Tadeusz  Kulczycki}
\address{Institute of Mathematics, Wroc{\l}aw
University of Technology, Wyb. Wyspianskiego 27, 50-370 Wroc{\l}aw, Poland}
\email{tkulczyc@im.pwr.wroc.pl}
\thanks{T. Kulczycki was supported by KBN
grant 2 P03A 041 22
and RTN Harmonic Analysis and  Related Problems, contract
HPRN-CT-2001-00273-HARP}
\author{Pedro J.
M\'endez-Hern\'andez}
\address{ Department of Mathematics, The University of Utah, 155 S. 1400 E. Salt lake City,
UT, 84112--0090}
\email{mendez@math.utah.edu}
 \keywords{Symmetric stable processes,
ground state eigenfunctions, multiple integrals} \subjclass{30C45}
\begin{abstract}
{\it We prove that the ground state eigenfunction for
symmetric stable processes of order $\alpha\in (0, 2)$ killed upon
leaving the interval $(-1, 1)$  is concave on $(-\frac{1}{2},
\frac{1}{2})$.  We call this property  ``mid--concavity."  A
similar statement holds for rectangles in $\R^d$, $d>1$. These
result follow  from similar results for finite dimensional
distributions of Brownian motion and subordination. }
\end{abstract}

\maketitle

\section{Introduction}

Let $D$ be a bounded convex domain in  $\R^d$, $d\geq 1$, and let
$\varphi_1$ be the first eigenfunction for the Dirichlet Laplacian
in $D$. In their seminal paper \cite {BrLi}, Brascamp and Lieb
proved that $\varphi_1$  is log--concave in $D$. That is,
$\log(\varphi_1)$ is concave on any segment contained in the
domain. This result has led to many interesting applications in
analysis, geometry, pde, mathematical physics and probability. For
some of these applications, see  Borell \cite{Bor1}, \cite{Bor2},
\cite{Bor3} and the many references therein.  In particular, the
log--concavity of $\varphi_1$ leads to  estimates of the spectral
gap  $\lambda_2-\lambda_1$ which in tern describe the rate to
equilibrium of the Brownian motion conditioned to remain forever
in the domain $D$. We refer the reader to \cite{BM}, \cite{Ling},
\cite{singer} and \cite{smits} for some of these applications and
additional references.

In \cite{BaKu}, the first two authors initiated the study of what may be
 called the
``fine spectral theoretic properties" of symmetric stable
processes.  Unfortunately, given the ``nonlocality"  of the
generator of these processes, even the most basic questions seem
to be very difficult.  It was  proved
%tkchange
 in \cite{BaKu} (Theorem 5.1) that the  ground state  eigenfunction
for  the Cauchy
process in the interval $(-1, 1)$
is concave. We, of course, expect  this to be the case for any
symmetric stable process.  The purpose of this paper is to prove
that for any symmetric stable processes, the ground state
eigenfunction is concave in $(-\frac{1}{2}, \frac{1}{2})$.  We
call this property  {\it ``mid--concavity".}  This will follow from a more general result
on ``mid--concavity" of the finite dimensional distributions of these processes.
This ``mid--concavity" result is new even for Brownian motion.

We first recall some basic definitions. Let $X_t^{\alpha}$ be a d-dimensional symmetric stable
process of index $0<\alpha \leq 2$.
The process $X_t^{\alpha}$ has stationary independent increments and its
transition density $p^{\alpha}_t(x,y)=p^{\alpha}_t(x-y)$, $t > 0$, $x,y \in
\R^d$, is
determined by its Fourier transform
$$
\exp(-t|\xi|^\alpha) = \int_{\R^d} e^{i \xi\cdot
y}p^{\alpha}_t(y)\,dy.
$$
 These are L\'evy processes with  right
continuous sample paths.  The  transition densities satisfy the
scaling  property $$ p^{\alpha}_t(x, \, y)= t^{-d/\alpha}
p^{\alpha}_1(t^{-1/\alpha}x,\,t^{-1/\alpha} y), $$ hence the
process has the scaling property of index $\alpha$. When
$\alpha=2$, $X_t^2$ is just Brownian motion $B_t$ running at twice
the speed and when $\alpha=1$, $X_t^1$ is the Cauchy process.   In
the first case,  $p^2_t(x, y)$ is the usual Gaussian distribution
(heat kernel)  and in the second, $ p^1_t(x,y)$ is the Cauchy
distribution (Poisson kernel).

Our interest here is on symmetric stable   processes of index
$0<\alpha<2$ killed  upon leaving a domain $D$. That is, let $D
\subset \Rd$, $d\geq 1$, be a  nonempty bounded connected open set
and let $$\tau_{D}^{\alpha} = \inf\{t \ge 0: X_{t}^{\alpha} \notin
D\}$$ be the first exit time of $X^{\alpha}_t$ from  $D$. Let
$$T_{t}^{D}f(x) = E_{x}(f(X_{t}^{\alpha}), \tau^{\alpha}_{D} > t),
$$ for $x\in D$, $t>0$ and $f\in L^{2}(D)$,  be the semigroup of
the killed process. The  killed process has  transition densities
$p_{D}^{\alpha}(t,x,y)$ and
\begin{equation}
T_{t}^{D}f(x) = \int_{D} p_{D}^{\alpha}(t,x,y) f(y) \, dy.
\end{equation}
 As with Brownian motion,
\begin{equation}
p_{D}^{\al}(t, x,y) = p^{\alpha}(t, x,y) - r_{D}(t, x,y),
\end{equation}
where
\begin{equation}
r_{D}(t, x,y) = E_{x}(p^{\al}_{t - \tau_{D}^{\alpha}}( X^{\alpha}_{\tau_{D}^{\alpha}},y),
\tau^{\alpha}_{D} < t).
\end{equation}
 From this it follows  that the transition  function $p_{D}^{\al}(t, x,y)$   is nonnegative,
symmetric, jointly continuous in $x$ and $y$, and that for all
$x, y\in D$ and
$ t>0$,
%tkchange
%$p_D^{\al} (t, x, y)\leq p^{\alpha}_t(x, y)=t^{-d/\alpha}p^{\al}_1(x/t,y/t)
%\leq Ct^{-d/\alpha},
%$
$$p_D^{\al} (t, x, y)\leq p^{\alpha}_t(x, y)=t^{-d/\alpha}
p^{\al}_1(t^{-1/\alpha} x,t^{-1/\alpha} y)
\leq C t^{-d/\alpha},
$$
where $C=(2\pi)^{-d}\omega_{d}\G(d/\al)/\al$ and $\omega_{d}$ is the
surface measure of the unit sphere in $\Rd$. In fact,  $p_{D}^{\alpha}(t,x,y)$ is strictly positive for
$x, y\in  D$.  These properties and
the  general theory of  heat semigroups
(as in \cite{Da1}) gives an
orthonormal basis of eigenfunctions
$\{\varphi_n^{\alpha}\}$
  on
$L^2(D)$ with  eigenvalues $\{\lambda_n^{\alpha}\}$ satisfying
$0<\lambda_1^{\alpha}<\lambda_2^{\al}\leq \lambda_3^{\al}\leq
\dots,$ and $\lambda_n^{\al}\to\infty$, as $n\to\infty$. That is,
$$T_t^{D}{\varphi_n^{\alpha}}(x)=e^{-\lambda_n^{\al}t}\varphi_n^{\alpha}(x), \, x\in D.
$$
 In addition, the first eigenvalue $\lambda_1^{\al}$ is simple
and its corresponding eigenfunction $\varphi_1^{\al}$, which we
will refer to as the {\it ground  state eigenfunction},  is an
analytic strictly positive function on $D$.  The infinitesimal
generator of the
%tkchange
%semigroup is $(-\Delta)^{\al/2}$. We can think of the
%eigenfunction and eigenvalues as solutions to the eigenvalue
%problem $(-\Delta)^{\alpha/2} \varphi_{n}(x)=\lambda_{n}^{ \al}
%\varphi_{n}^{\alpha}(x)$, $x\in D$ and $\varphi_{n}^{\alpha}(x)
%=0$ for $x\in D^c$; the Dirichlet problem for stable processes.
%We refer the reader to \cite{BG1}, \cite{BG2} and \cite{G1} where
%many of the general properties of the $\alpha$--stable semigroup
%are established.
semigroup is $-(-\Delta)^{\al/2}$. We can think of the
eigenfunction and eigenvalues as solutions to the eigenvalue
problem
$$(-\Delta)^{\alpha/2} \varphi_{n}^{\al}(x)=\lambda_{n}^{ \al}
\varphi_{n}^{\alpha}(x),
$$
 $x\in D$ and $\varphi_{n}^{\alpha}(x)
=0
$
for $x\in D^c$; the Dirichlet problem for stable processes. We
refer the reader to \cite{BG1}, \cite{BG2}, \cite{BB}, \cite{CS}
and  \cite{G1} where many of the general properties of the
$\alpha$--stable semigroup and its generator are established.

 The following question is motivated from  the
result of Brascamp and Lieb \cite{BrLi} mentioned above for
Brownian motion and by  its many applications.

\begin{question}\label{question1}
Let $D\subset \R^d$, $d\geq 1$,  be a bounded convex domain and
 $0<\alpha<2$. Is $\varphi_1^{\alpha}$
log--concave? In other words,   is $\log(\varphi_1^\alpha)$  concave on any
segment contained in $D$?
\end{question}

The only  known case   is  when    $D=(-1, 1)$ and $\alpha=1$,
where the question is  answered in the affirmative in \cite{BaKu}.
In fact, it is shown  in \cite{BaKu} that the ground state
eigenfunction for the Cauchy process in $(-1, 1)$ is concave.
Because  of this case we believe this result should hold for all
$\alpha$--stable processes.  More precisely, we have

\begin{conjecture}\label{conjecture1} Let $\varphi_1^\alpha$ be the ground state
  eigenfunction for the
symmetric stable processes of index $0<\alpha<2$ killed upon
leaving  the interval $I=(-1, 1)$. Then $\varphi_1^\alpha$ is
concave on $I$.
\end{conjecture}

 There are by now many proofs of the log--concavity  result for  Brownian motion. None
of them, as far as we can see, adapt to the case of general
symmetric stable processes. However,  Brascamp--Lieb's  proof does
suggest some related questions which may provide some insight. We
briefly recall here their argument based on multiple integrals.
Let $B_t$ be Brownian motion and let $\tau_D$ be its first exit
time from
 $D$. Then one can show, see \cite{B}, that
$\varphi_1^2(x)=\lim_{t\to\infty}e^{\lambda_1^2 t}P_x\{\tau_D > 2
t\}$, uniformly in $x\in D$.  From this it is enough to prove that
$P_x\{\tau_D>t\}$ is log--concave in $x$ for every fixed $t>0$.
The latter can
  be written as the limit as $n$ and $k$ tend  to infinity  of
$P_x\{B_{jt/n}\in D_k; j=1, 2, \dots,  n\}$ where $D_k$ is a
sequence of convex domains strictly increasing ($\overline
D_k\subset D_{k+1}$) up to $D$. We then  reduce the problem to
prove that for any
 convex domain $D$, $P_x\{B_{jt/n}\in D;  j=1, 2, \dots, n\}$ is
log--concave on $D$ as a function of $x$,  for all $t>0$ and all $n$.
 This, however, is a multiple convolution of Gaussians with the indicator function of the set
$D$.  Since the Gaussian $p_t^2(x)$ is log--concave for all $t>0$ and the indicator function
 of a convex domain is log--concave, the result follows
from the fact that convolutions of log--concave functions are log--concave.
 Using right continuity of paths, we can try to  repeat  this argument  for $\alpha$--stables
processes. However, this time the argument breaks down right at the end.
 For example, if $\alpha=1$ the density for the
 Cauchy process, $p_t^1(x, y)=p_t^1(x-y),$
is not log--concave for all $t$.  The obvious  variation of this argument using the fact that
$X_t^{\alpha} = B_{2\sigma_t}$, where $\sigma_t$ is a  stable subordinator of index $\alpha /2$
independent of $B_t$,   also fails basically due to the fact that  the sum of
log--concave functions is not necessarily log--concave.

There is however, a substitute for log--concavity which gives some  insight into
the shape of the ground state eigenfunction.
%tkchange
%for symmetric convex
%domains in the plane and
%for the purpose of obtaining some
%spectral gap estimates,  suffices.
We call this property \emph{``mid--concavity".}

\begin{definition}\label{definition1}
Let $D\subset \R^d$ be a  convex domain which is symmetric
relative to each coordinate  axes.  Let $J$ be a line  segment in
$D$ parallel to the $x_1$-axis which intersects the boundary
$\partial D$ only at the two points $(-a_1, a_2, \dots, a_d)$,
$(a_1, a_2, \dots, a_d)$, $a_1>0$.  We will say that the function  $F:D\to \R$,
 is  \emph{mid--concave} on $J$ if it is concave on the segment
(half of $J$) from the point $(-a_1/2, a_2, \dots, a_d)$ to
$(a_1/2, a_2, \dots, a_d)$.  The function is mid--concave along
the $x_1$--axis if it is mid--concave on every such segment
contained in $D$ which is parallel to the $x_1$--axis.  A similar
definition applies for {\it mid--concavity} along the $x_2$-axis,
$\cdots$, $x_d$--axis.
 The function is mid--concave on $D$ if it is mid--concave
along  each coordinate axes.
\end{definition}

Our main result in this paper is the following
\begin{theorem}\label{thm1}
Let $Q=(-a_1, a_1)\times (-a_2, a_2)\times \cdots \times(-a_d, a_d)$, $0<a_i<\infty$
for all $i=1, 2, \dots, d$,
be a rectangle in $\R^d$.  The ground state eigenfunction $\varphi_1^{\alpha}$
for the symmetric stable
process of index $0<\alpha<2$
 is {\it mid--concave} on $Q$. In addition, if
$x=(x_1,\ldots,x_n)  \in Q$, then
\begin{equation}\label{mono1}
\frac{\partial}{\partial x_i}  \varphi_1^{\alpha}(x)\geq 0, \text{
if } x_i<0, \text { and } \frac{\partial}{\partial x_i}
\varphi_1^{\alpha}(x)\leq 0, \text{ if } x_i>0. \end{equation}
\end{theorem}
 Using arguments of multiple integrals as described above,  we will show that Theorem
\ref{thm1} follows from

\begin{theorem}\label{thm2}
Let $Q$ be a rectangle in $\R^d$.  Let $0<t_1<t_2< \dots <t_n<\infty$. The function
\begin{equation}
F(x)=P_x\{X_{t_1}^{\alpha}\in Q, \dots, X_{t_n}^{\alpha}\in Q\}
\end{equation}
is mid--concave in $Q$ for any $0<\alpha\leq 2$.  In addition, if
$x=(x_1,\ldots,x_n)  \in Q$, then
\begin{equation}\label{mono2}
\frac{\partial}{\partial x_i}  F(x)\geq 0, \text{ if } x_i<0,
\text { and } \frac{\partial}{\partial x_i} F(x)\leq 0, \text{ if
} x_i>0. \end{equation}
\end{theorem}
\begin{remark}\label{remark1} It is important to note here that Theorem
 \ref{thm2} is new even in the Brownian motion case ($\alpha=2$).   Indeed, as we shall see, the case $\alpha=2$ implies
the general case by subordination.
\end{remark}
%tkchange
%For  any  $n=1, 2, \dots$, let
%$0<t_1<t_2<
%\dots <t_n<\infty$. With this line of reasoning  if we could prove that the
%fucntion
%\begin{equation}
%F(x)=P_x\left\{\,X_{t_1}^\alpha\in I, \dots, X_{t_n}^\alpha \in
%I\right\}
%\end{equation}
%is  convex on $I$,  the  Conjecture \ref{conjecture1}  would follow not only
%for the
%interval but also for any rectangle in $\R^d$, $d>1$.
%Unfortunately, such an $F$, while {\it mid--concave} on $I$, is not
%concave on all of $I$, as we will show below.

If we consider the eigenfunction for the Laplacian in the unit
 disk $\bD$ in the plane,  one can
show,  by analysis of the Bessel function, that such a function
is not concave in $\bD$\,  but it is
{\it mid--concave.} Also, it may be tempting to conjecture that for
 any symmetric domain in the plane
the eigenfunction is  {\it mid--concave.}
This, however,  is not the case, even for the Brownian motion,
 as we will show at the end of the paper.

The paper is organized as follows.  In $\S2$, we prove
 that the multiple convolutions of Gaussians in the interval
$(-1, 1)$ is {\it mid--concave}.  In $\S3$, we
 show how this and subordination implies Theorem \ref{thm2}.
 Here we also show that full concavity fails for general multiple integrals
 and that {\it mid--concavity} fails in general symmetric domains in the plane.

\section{Mid--concavity for Brownian motion}

Let \[p_t(x)=\frac{1}{\sqrt{2\pi t}}e^{-\frac{x^2}{2t}}\] be the
Gaussian density in one dimension. With the notation of the
introduction, we have $p^2_t (x, y)= p_{2t}(x-y)$.

\begin{proposition}\label{proposition1} Let $n=1, 2, \dots$ and let
$t_1, t_2, \dots, t_n$ be  real numbers in $(0, \infty)$. For $x\in (-1, 1)$  define
\begin{equation}\Phi_n(x)= \int_{-1}^{1}\cdots\int_{-1}^{1}\, \prod_{i=1}^n
p_{t_i}(x_{i-1}-x_i)\,dx_1\ldots dx_n,
\end{equation}
where $x_0=x$.   The function $\Phi_n(x)$ is {\it mid--concave}
on $(-1, 1)$.  That is, $\Phi_n(x)$ is concave on $(-\frac{1}{2},
\frac{1}{2})$.
\end{proposition}

Clearly   $\Phi_n(x)$ is a positive  even function on $[-1,1]$.
Integrating  by parts we obtain
\begin{equation}\begin{split}
-\frac{\partial}{\partial x} \Phi_n(x)&=\, \frac{1}{\sqrt{2\pi
t_n}}\, \int_{-1}^{1}\;\left(\frac{\partial}{\partial y}
e^{-\frac{(y-x)^2}{2t_n}}\right)\, \Phi_{n-1}(y)\, dy \\& =
\frac{\Phi_{n-1}(1)}{\sqrt{2\pi t_n}} \left(
e^{-\frac{(1-x)^2}{2t_n}}-e^{-\frac{(1+x)^2}{2t_n}}\right)\\&- \,
\frac{1}{\sqrt{2\pi t_n}}\,\int_{-1}^1\,
e^{-\frac{(y-x)^2}{2t_n}}\, \frac{\partial }{\partial
y}\Phi_{n-1}(y)\, dy .
\end{split}\label{parts2}\end{equation}
Notice that for all $t>0$,
\begin{equation} \left(
e^{-\frac{(1-x)^2}{2t}}-e^{-\frac{(1+x)^2}{2t}}\right)=\,e^{\frac{-(1-x)^2}{2}}\,\left(
1-e^{-\frac{2x}{t}}\right),\label{parts1}\end{equation}   is a
positive increasing  function on $[0,1]$.

\begin{lemma}\label{lemma1}
The function $\Phi_n(x)$ is decreasing on $(0,1)$ for all $n\geq
1$.
\end{lemma}
\begin{proof}
We argue by induction. If $n=1$, then
\begin{equation}\begin{split}
\frac{\partial}{\partial x} \Phi_1(x)&=\,\frac{1}{\sqrt{2\pi
t_1}}\, \int_{-1}^{1}\frac{-\partial}{\partial y}
e^{-\frac{(y-x)^2}{2t_1}}\,dy
\\& = \frac{1}{\sqrt{2\pi t_1}} \left(
-e^{-\frac{(1-x)^2}{2t_1}}+e^{-\frac{(1+x)^2}{2t_1}}\right)\\& <
0,
\end{split}\end{equation}
for all $ x \in (0,1)$. Thus $\Phi_1(x)$ is decreasing on  $(0,1)$.

Let us assume that $\Phi_{n-1}(x)$ is decreasing on $(0,1)$.  That is, suppose that
\[\frac{\partial }{\partial x}\Phi_{n-1}(x) \leq 0,\]
for all $x \in (0,1)$.
Because of (\ref{parts1}), it is enough to prove that
\begin{equation} \, \frac{-1}{\sqrt{2\pi t_n}}\,\int_{-1}^1
e^{-\frac{(y-x)^2}{2t_n}} \frac{\partial }{\partial
y}\Phi_{n-1}(y) \,dy \geq 0. \label{decreasing}\end{equation} By
symmetry  \[\frac{\partial }{\partial
y}\Phi_{n-1}(y)=-\frac{\partial }{\partial y}\Phi_{n-1}(-y).\] On
the other hand,  if $x>0$ then
\[ e^{ -\frac{(x-y)^2}{2t}} \geq e^{-\frac{(x+y)^2}{2t}},\]
for all $t, y >0$. Hence for all $y>0,$
\[-e^{-\frac{(x-y)^2}{2t_n}}\frac{\partial }{\partial
y}\Phi_{n-1}(y)- e^{-\frac{(x+y)^2}{2t_n}}\frac{\partial
}{\partial y}\Phi_{n-1}(-y) =\]\[
-e^{-\frac{(x-y)^2}{2t_n}}\frac{\partial }{\partial
y}\Phi_{n-1}(y)+ e^{-\frac{(x+y)^2}{2t_n}} \frac{\partial
}{\partial y}\Phi_{n-1}(y) = \] \[\,\left(
-e^{-\frac{(x-y)^2}{2t_n}}+ e^{-\frac{(x+y)^2}{2t_n}}\,\right)\,
\frac{\partial }{\partial y}\Phi_{n-1}(y) \geq 0 .\] Integrating
this inequality on $[0,1]$ we obtain (\ref{decreasing}).
\end{proof}

Notice that
\begin{equation}\begin{split}
\frac{\partial^2}{\partial x^2} \Phi_1(x)&=\, \frac{1}{\sqrt{2\pi
t_1}}\, \int_{-1}^{1}\frac{\partial^2}{\partial y^2}
e^{-\frac{(y-x)^2}{2t_1}}\,dy \\& = \frac{1}{t_1\sqrt{2\pi t_1}}
\left( -(1-x)
e^{-\frac{(1-x)^2}{2t_1}}-(1+x)e^{-\frac{(1+x)^2}{2t_1}}\right)\\&
< 0, \label{n=1}\end{split}\end{equation} for all $ x \in (-1,1)$.
Thus $\Phi_1(x)$ is concave in $(-1,1)$. We will know prove that
$\Phi_n(x)$ is concave in $(-\frac{1}{2},\frac{1}{2})$.

\begin{lemma}\label{lemma2} If $0 \leq x \leq \frac{1}{2}$, then  for all $n \geq 1,$ \[
\frac{\partial}{\partial x}\Phi_n(x) \geq \frac{\partial}{\partial
x}\Phi_n(1-x).\]
\end{lemma}
\begin{proof}
By (\ref{n=1}) the result is true for $n=1$. Let us assume that
the result is true for $n-1$.  Let
\[ \psi_n(x)= \,
\frac{-1}{\sqrt{2\pi t_n}}\,\int_{-1}^1 e^{-\frac{(y-x)^2}{2t_n}}
\frac{\partial }{\partial y}\Phi_{n-1}(y)\, dy .\] Because of
(\ref{parts1}), it is enough to prove that

\begin{equation}\label{psi}
 \psi_n(1-x)\geq \psi_n(x).\end{equation}
Let $ y \in (-1,0)$, then \[1-x-y\geq x-y\geq 0 .\] Thus
\[e^{-\frac{(1-x-y)^2}{2t_n}}\leq e^{-\frac{(x-y)^2}{2t_n}}.\]
Since
\[-\frac{\partial }{\partial y}\Phi_{n-1}(y) \leq 0,\]
for all $y < 0$,   it follows that
\[-\int_{-1}^0 e^{-\frac{(x-y)^2}{2t_n}}
\frac{\partial }{\partial y}\Phi_{n-1}(y) \,dy \leq -\int_{-1}^0
e^{-\frac{(1-x-y)^2}{2t_n}} \frac{\partial }{\partial
y}\Phi_{n-1}(y)\, dy.
\]
To simplify notation let \[ \phi(y)= -\frac{\partial }{\partial
y}\Phi_{n-1}(y).\] Let $ y \in (0,\frac{1}{2})$, and consider
$\hat{y}=1-y$. Notice that $ \hat{y} \in (\frac{1}{2},1)$ and \[
\frac{1}{2}-y=\hat{y}-\frac{1}{2}.\] By induction,
\[0\leq \phi(y) \leq \phi(\hat{y}). \]
On the other hand,
\begin{eqnarray*}
e^{-\frac{(x-y)^2}{2t_n}}&=e^{-\frac{(\hat{x}-\hat{y})^2}{2t_n}},\\
e^{-\frac{(x-\hat{y})^2}{2t_n}}&=e^{-\frac{(\hat{x}-y)^2}{2t_n}},\\
e^{-\frac{(x-y)^2}{2t_n}}&\geq
e^{-\frac{(x-\hat{y})^2}{2t_n}}.\end{eqnarray*} Thus
\[ \left( e^{-\frac{(x-y)^2}{2t_n}}-e^{-\frac{(\hat{x}-y)^2}{2t_n}}\right) \,\phi(y)
\,\leq\,\left(
e^{-\frac{(\hat{x}-\hat{y})^2}{2t_n}}-e^{-\frac{(x-\hat{y})^2}{2t_n}}\right)
\phi(\hat{y}),\] and we conclude that
\[  e^{-\frac{(x-y)^2}{2t_n}}\phi(y)+ e^{-\frac{(x-\hat{y})^2}{2t_n}}
\phi(\hat{y}) \leq e^{-\frac{(\hat{x}-y)^2}{2t_n}} \phi(y)+
e^{-\frac{(\hat{x}-\hat{y})^2}{2t_n}} \phi(\hat{y}).\] Integrating
over $(0,\frac{1}{2})$ we obtained that
\[-\int_{0}^1 e^{-\frac{(x-y)^2}{2t_n}}
\frac{\partial }{\partial y}\Phi_{n-1}(y)\, dy \leq -\int_{0}^1
e^{-\frac{(1-x-y)^2}{2t_n}} \frac{\partial }{\partial
y}\Phi_{n-1}(y)\, dy,
\]the desired result immediately follows.
\end{proof}

\begin{lemma}\label{lemma3} If $0 \leq x < u \leq \frac{1}{2}$, then  for all $n \geq 1,$ \[
\frac{\partial}{\partial x}\Phi_n(x) \geq \frac{\partial}{\partial
x}\Phi_n(u).\]
\end{lemma}
\begin{proof}
By (\ref{n=1}) the result is true for $n=1$. Let us assume that
the result is true for $n-1$.  As in Lemma \ref{lemma1},  we let
\[ \psi_n(x)= \,
\frac{-1}{\sqrt{2\pi t_n}}\,\int_{-1}^1 e^{-\frac{(y-x)^2}{2t_n}}
\frac{\partial }{\partial y}\Phi_{n-1}(y) \,dy \] and
\[ \phi(y)= -\frac{\partial }{\partial
y}\Phi_{n-1}(y).\]
By  (\ref{parts1}), it is enough to prove that
\begin{equation}\label{psi2}
\psi_n(u)\geq \psi_n(x).
\end{equation}
Let $ y \in (-1,0)$, then $|u-y|\geq |x-y|$. Thus
\[e^{-\frac{(u-y)^2}{2t_n}}\leq e^{-\frac{(x-y)^2}{2t_n}}.\]
Lemma \ref{lemma1} implies that
\[-\int_{-1}^0 e^{-\frac{(x-y)^2}{2t_n}}
\frac{\partial }{\partial y}\Phi_{n-1}(y)\, dy \leq -\int_{-1}^0
e^{-\frac{(u-y)^2}{2t_n}} \frac{\partial }{\partial
y}\Phi_{n-1}(y) \,dy.
\]

Let $m=\frac{x+u}{2}$. For all  $ y \in (0,m)$, define
$\tilde{y}=x+u-y$. Notice that $\tilde{y} \in (m,x+u)$ and that
\[ |x-y|=|u-\tilde{y}|.\]

We can  easily check that
\begin{eqnarray*}
e^{-\frac{(x-y)^2}{2t_n}}&=e^{-\frac{(u-\tilde{y})^2}{2t_n}},\\
e^{-\frac{(x-\tilde{y})^2}{2t_n}}&=e^{-\frac{(u-y)^2}{2t_n}},\\
e^{-\frac{(x-y)^2}{2t_n}}&\geq
e^{-\frac{(x-\tilde{y})^2}{2t_n}},\end{eqnarray*} for all $y \in
(0,m)$. We claim that

\begin{equation}\phi(y) \leq \phi(\tilde{y}).
\end{equation}
This follows immediately from the induction hypothesis if
$\tilde{y}\leq \frac{1}{2}$. On the other hand, if
$\tilde{y}=(x+u)-y\geq \frac{1}{2}$, then
\[1-(x+u)+y \leq \frac{1}{2},\; \text{and}\;y \leq 1-(x+u)+y.\]  Lemma \ref{lemma2} and the
induction hypothesis imply  that

\[0\leq \phi(y) \leq \phi( 1-(x+u)+y) \leq \phi( (x+u)-y)= \phi(\tilde{y}). \]
Thus

\[  e^{-\frac{(x-y)^2}{2t_n}}\phi(y)+ e^{-\frac{(x-\tilde{y})^2}{2t_n}}
\phi(\tilde{y}) \leq e^{-\frac{(u-y)^2}{2t_n}} \phi(y)+
e^{-\frac{(u-\tilde{y})^2}{2t_n}} \phi(\tilde{y}).\] Integrating
over $(0,m)$ we obtained that

\[-\int_{0}^{x+u} e^{-\frac{(x-y)^2}{2t_n}}
\frac{\partial }{\partial y}\Phi_{n-1}(y)\, dy \leq
-\int_{0}^{x+u} e^{-\frac{(u-y)^2}{2t_n}} \frac{\partial
}{\partial y}\Phi_{n-1}(y)\, dy.
\]

Finally if $ y\in [x+u,1]$ then
\[ e^{-\frac{(x-y)^2}{2t_n}}\leq e^{-\frac{(u-y)^2}{2t_n}}.\]
Therefore
\[-\int_{x+u}^{1} e^{-\frac{(x-y)^2}{2t_n}}
\frac{\partial }{\partial y}\Phi_{n-1}(y) \,dy \leq
-\int_{x+u}^{1} e^{-\frac{(u-y)^2}{2t_n}} \frac{\partial
}{\partial y}\Phi_{n-1}(y) \,dy.
\]
\end{proof}
By symmetry, Proposition \ref{proposition1} follows from Lemma
\ref{lemma3}. The following is an immediate corollary of
Proposition \ref{proposition1}.

\begin{corollary}\label{corollary1} Let $B_t$ be one dimensional Brownian motion and set
$I=(-a, a)$,  For
$0<t_1<t_2<\dots <t_n$,  the function
\begin{equation}
F(x)=P_x\{B_{t_1}\in I, B_{t_2}\in I, \dots, B_{t_n}\in I\}
\end{equation}
is mid--concave in $I$. In addition, if $x \in I$, then
\begin{equation} \label{mono3}F'(x)\geq 0, \text{ if } x<0, \text { and }
F'(x)\leq 0, \text{ if } x>0.
\end{equation}
\end{corollary}

\begin{proof}
By the Markov property,
\begin{equation}
F(x)=\int_{-a}^{a}\cdots\int_{-a}^{a}\, \prod_{i=1}^n
p_{t_i-t_{i-1}}(x_{i-1}-x_i)\,dx_1\ldots dx_n,
\end{equation}
where $x_0=x$ and $t_0=0$.  This is exactly the same expression as
in Lemma 2.1 and Proposition \ref{proposition1} except for the
fact that the interval $(-1, 1)$ has been replaced by  the
interval $(-a, a)$. The proof of the proposition is the same for
this case and the corollary follows.
\end{proof}

\begin{corollary}\label{corollary2}
Let $B_t$ be Brownian motion in $\R^d$ and let
$Q=I_1\times I_2\times\dots \times I_d$ where $I_i=(-a_i, a_i),$ be a rectangle in $\R^d$.
 For
$0<t_1<t_2<\dots <t_n$,  the function
\begin{equation}
F(x)=P_x\{B_{t_1}\in Q, B_{t_2}\in Q, \dots, B_{t_n}\in Q\}
\end{equation}
is mid--concave in $Q$. In addition, if  $x=(x_1, x_2, \dots, x_d)
\in Q$, then
\begin{equation} \label{mono4}
\frac{\partial}{\partial x_i}F(x)\geq 0, \text{ if } x_i<0, \text
{ and } \frac{\partial}{\partial x_i} F(x)\leq 0, \text{ if } x>0.
\end{equation}

\end{corollary}
\begin{proof}
With $x=(x_1, x_2, \dots, x_d)$ and $B_t=(B_t^1, B_t^2, \dots,
B_t^d)$, it follows by independence that \[F(x)=\prod_{i=1}^d
P_{x_i}\{B^i_{t_1}\in I_i, B^i_{t_2}\in I_i, \dots, B^i_{t_n}\in
I_i\}\] and the conclusion of the corollary follows from Corollary
\ref{corollary1} and our definition of {\it mid--concavity} for
domains in $\R^d$.
\end{proof}

\section{Mid--concavity for stable processes}

In this section we prove Theorems \ref{thm1} and  \ref{thm2}.
First, let us recall that for  $ 0< \alpha
<2$ the symmetric stable process $X_t^{\alpha}$ in $\R^d$ has the representation
\begin{equation}X_t^{\alpha} = B_{2\sigma_t}\label{3equa1},
\end{equation}
where $\sigma_t$ is a  stable subordinator of index $\alpha /2$
independent of $B_t$ (see \cite{BG3}). Thus
\begin{equation}\label{subordination}
p_t^{\alpha}(x-y)= \int_0^{\infty} p_s^2(x-y)
g_{\alpha/2}(t,s)ds,
\end{equation} where $ g_{\alpha/2}(t,s)$ is the
transition density of $\sigma_t$ and
\[p^2_{t/2}(x-y)=\frac{1}{({2\pi t})^{d/2}}e^{-\frac{|x-y|^2}{2t}}.\]
%tkchange
%As before, we will simply write $p_t(x-y)$ for $p_{t}^2(x-y)$.

%Now, let $Q$, $\alpha$ and $t_1, t_2, \dots, t_n$ be as
Now, let $Q$ and $t_1, t_2, \dots, t_n$ be as in the statement of
Theorem \ref{thm2}.  Set $x_0=x$ and $t_0=0$.  Using the Markov
property of the stable process $X_t^{\alpha}$, the subordination
formula (\ref{subordination}), Fubini's theorem, and the Markov
property of the Brownian motion, in this order, we obtain,
\begin{eqnarray*}
&& F(x)=P_x\{X_{t_1}^{\alpha}\in Q, \dots, X_{t_n}^{\alpha}\in
Q\}\\ && = \int_Q\cdots\int_Q\, \prod_{i=1}^n
p^{\alpha}_{t_i-t_{i-1}}(x_{i-1}-x_{i})\,dx_1\ldots
dx_n\nonumber\\ && =\int_0^{\infty}\ldots\int_0^{\infty}\left(
\int_Q\cdots\int_Q\, \prod_{i=1}^n
p_{s_i}^2(x_{i-1}-x_{i})\,dx_1\ldots dx_n\right)\\ &&\times
\prod_{i=1}^{n}g_{\alpha/2}(t_i-t_{i-1},s_i)\,ds_1\ldots ds_n\\
&& = \int_0^{\infty}\ldots\int_0^{\infty}P_x\{B_{2 s_1}\in Q,
B_{2(s_1+s_2)}\in Q, \dots, B_{2(s_1+s_2+\dots +s_n)}\in Q\}\\ &&\times
\prod_{i=1}^{n}g_{\alpha/2}(t_i-t_{i-1},s_i)\,ds_1\ldots ds_n.
\end{eqnarray*}
Since the function $$ P_x\{B_{2s_1}\in Q, B_{2(s_1+s_2)}\in Q,
\dots, B_{2(s_1+s_2+\dots +s_n)}\in Q\} $$ is {\it mid--concave}
and satisfies the monotonicity property (\ref{mono4}), by
Corollary \ref{corollary2}, so is the integral against the
densities $g_{\alpha/2}(t_i-t_{i-1},s_i)$ and this completes the
proof of Theorem \ref{thm2}.

With Theorem \ref{thm2} proved, we argue as in the proof of the
{\it log--concavity} for Brownian motion discussed in the
introduction.  Recall that $\varphi_1^{\alpha}$ is the ground
state eigenfunction for the stable process of index $\alpha$, $
\alpha \in (0, 2)$,  killed upon leaving $Q$ and
$\lambda_1^{\alpha}$ is its eigenvalue. Let $\tau_Q^{\alpha}$ be
the first exit time of the symmetric stable process from $Q$.
Since $Q$ is certainly intrinsically ultracontractive, see
\cite{Ku}, we have that
\begin{equation}
\varphi_1^{\alpha}(x)=\lim_{t\to\infty}e^{\lambda_1^{\alpha} t}P_x\{\tau_Q^{\alpha}>t\}.
\end{equation}
The convergence is uniform for $x\in Q$.   Thus to prove {\it
mid--concavity} for $\varphi_1^{\alpha}(x)$ it is enough to prove
mid--concavity for $P_x\{\tau_Q^{\alpha}>t\}$. By the right
continuity of the sample paths, we have,
\begin{eqnarray}\label{approx}
P_x\{ \,\tau_Q^{\alpha} >t\,\} &=& P_z\{\,  X_{s}^{\alpha} \in Q,0\leq s
\leq t\,\}\nonumber\\
&=&\lim_{n \to \infty} P_x \{\, X_{\frac{it}{n}}^{\alpha} \in
Q, i= 1,\ldots,n\,\}.
\end{eqnarray}
Theorem \ref{thm1} now follows from this and Theorem \ref{thm2}.

We remark that in the case of Brownian motion,  there is an extra approximation by an
increasing sequence of domains in passing from the first equality to the second in
(\ref{approx}). This is not needed for our stable processes since,  as
explain in
%tkchange
%\cite{BaKu}, for any domain
\cite{Bo}, Lemma 6, for any domain $D\subset\R^d$ with Lipschitz
boundary, \[P_x\{X^{\alpha}_{\tau_D}\in \partial D\}=0\, \text{
for } x\in D.\]

The above argument applies not only to symmetric stable processes but also to any other
process which is obtained by subordination of Brownian motion.  In particular, the above
results hold for the so called  ``relativistic" Brownian motion and ``relativistic"
$\alpha$--stable processes studied in
\cite{ry}.

It is of course natural to ask if the function of Proposition
\ref{proposition1} is concave in the whole interval $(-1, 1)$ for
all $n$ and all $t_i$. Notice that, thanks to the proof of  Lemma
\ref{lemma2}, this is the case for $n=1$.
 If this were the case, it would show that the  same is true for the function
$P_x\{\,\tau_Q^{\alpha} >t\,\}$ and hence for the function
$\varphi_1^\alpha$, as desired. Unfortunately, this is not the
case.
\begin{proposition}\label{proposition2} Let
\begin{equation}\Phi_n(x)= \int_{-1}^{1}\cdots\int_{-1}^{1}\, \prod_{i=1}^n
p_{t_i}(x_{i-1}-x_i)\,dx_1\ldots dx_n,
\end{equation}
where $x_0=x$.  Then there exist a positive integer $n$ and  real
numbers $t_1, t_2, \dots, t_n$  in $(0, \infty)$ such that  the
function $\Phi_n(x)$ is {\it  not concave} on $(-1, 1)$.
\end{proposition}

\begin{proof}  We may replace, to simplify certain notation below,
the interval $(-1, 1)$ by the interval  $(0, \pi)$.  Fix $t$ and $s$ both positive.
Let $t_1=t$ and $t_2=\dots =t_n={\frac{s}{n-1}}$.  If the function $\Phi_n(x)$ is concave on
$(0, \pi)$ for all $n$  with these chosen $t_1, t_2, \dots, t_n$, letting $n\to\infty$ we see that the
function
\begin{eqnarray}
\int_0^{\pi}p_{t}(x-y)P_y\{\,\tau_{(0, \pi)} >s\,\}\,dy
\end{eqnarray}
is also concave on $(0, \pi)$.  Here we have used $\tau_{(0,
\pi)}$ to denote the first exit time of Brownian motion from the
interval.
%Since
%\begin{eqnarray}
%\lim_{s\to\infty}e^{\lambda_1 s}P_y\{\,\tau_{(0, \pi)} >s\,\}=
%\sqrt{\frac{2}{\pi}}\sin(y),
%\end{eqnarray}
%uniformly for $y\in (0, \pi)$, we see that for each $t>0$, the function
%\begin{eqnarray}
%F_t(x)=\int_0^{\pi}p_{t}(x-y)\sin(y)\,dy
%\end{eqnarray}
%must also be concave on $(0, \pi)$.
We have
\begin{eqnarray}
\lim_{s\to\infty}e^{\lambda_1 s}P_y\{\,\tau_{(0, \pi)} >s\,\}=
c \sin(y),
\end{eqnarray}
uniformly for $y\in (0, \pi)$, where $c > 0$ and $\lambda_1 = 1$ (the
first eigenvalue for $(0,\pi)$). It follows that for each $t>0$, the function
\begin{eqnarray}
F_t(x)=\int_0^{\pi}p_{t}(x-y)\sin(y)\,dy
\end{eqnarray}
must also be concave on $(0, \pi)$.

We will now show that the function $F_t(x)$ is not concave.  Without any difficulty we may differentiate under
the integral to obtain that
\begin{eqnarray}
F_t''(x)={\frac{1}{\sqrt{2\pi
}t^{5/2}}}\int_0^{\pi}\left[\,(x-y)^2-t\,\right]\,
e^{\frac{-(x-y)^2}{2t}} \sin(y)\,dy.
\end{eqnarray}
Taking $x=0$ and using the elementary inequality
\[y-{\frac{y^3}{3!}}\leq \sin(y)\leq y\] valid for all $y>0$, we
%tkchange
%see that
see that $F_t''(0)$ is equal to
\begin{eqnarray}
&&{\frac{1}{\sqrt{2\pi }t^{5/2}}}\,\int_0^{\pi}\left(y^2-t\right) \, e^{-\frac{y^2}{2t}} \sin(y)\,dy\nonumber\\
&=&{\frac{1}{\sqrt{2\pi }t^{5/2}}}\,\left(\int_0^{\pi}y^2\,
e^{-\frac{y^2}{2t}} \sin(y)\,dy-t\int_0^{\pi} e^{-\frac{y^2}{2t}}
\sin(y)\,dy\right)\nonumber\\ &\geq& {\frac{1}{\sqrt{2\pi
}t^{5/2}}}\left(\int_0^{\pi}y^3\,  e^{-\frac{y^2}{2t}}
\,dy-{\frac{1}{3!}}\int_0^{\pi}y^5\,  e^{-\frac{y^2}{2t}} \,dy-
t\int_0^{\pi}y e^{-\frac{y^2}{2t}}\, dy\right)\nonumber\\
&=&{\frac{1}{\sqrt{2\pi
}t^{5/2}}}\left(t^2\int_0^{\frac{\pi}{\sqrt{t}}}y^3\,
e^{-\frac{y^2}{2}}\,
dy-{\frac{t^3}{3!}}\int_0^{\frac{\pi}{\sqrt{t}}}y^5\,
e^{-\frac{y^2}{2}} \,dy- t^2\int_0^{\frac{\pi}{\sqrt{t}}}y
e^{-\frac{y^2}{2}} \,dy\right)\nonumber\\ &=& {\frac{1}{\sqrt{2\pi
t}}}\left(\int_0^{\frac{\pi}{\sqrt{t}}}y^3\, e^{-\frac{y^2}{2}}\,
dy-{\frac{t}{3!}}\int_0^{\frac{\pi}{\sqrt{t}}}y^5\,
e^{-\frac{y^2}{2}}\, dy- \int_0^{\frac{\pi}{\sqrt{t}}}y\,
e^{-\frac{y^2}{2}} \,dy\right).\nonumber
\end{eqnarray}

Since $$ {\frac{t}{3!}}\int_0^{\frac{\pi}{\sqrt{t}}}y^5\,
e^{-\frac{y^2}{2}} \,dy\to 0, $$

$$ \int_0^{\frac{\pi}{\sqrt{t}}}y^3\,  e^{-\frac{y^2}{2}}\,dy \to
2, $$ and $$ \int_0^{\frac{\pi}{\sqrt{t}}}y\,
e^{-\frac{y^2}{2}}\,dy \to 1, $$ as $t\to 0^+$, we see that
$F_t''(0)$ is positive for sufficiently small $t$.
%tkchange
%(In fact,
%$F_t''(0)=\infty$.)
By continuity, we have that $F_t''(x)>0$ for
sufficiently small $x\in (0, \pi)$ and sufficiently small $t$.
This, of course, contradicts the concavity of the function and
shows that $\Phi_{n}(x)$ is not concave.

\end{proof}

Of course, it may still be the case that the function $\Phi_n(x)$
is concave on the whole interval when we restrict to a sequence of
%tkchange
%times satisfying $t=t_1=t_2=\dots=t_n$ and substitute $p_t(x)$ by
%$p_t^\alpha(x)$, which what  is needed for   our
%applications.
times satisfying $t_1=t_2=\dots=t_n$ and substitute $p_{t_i}(x)$
by $p_{t_i}^\alpha(x)$, which is what  is needed for   our
applications (Conjecture 1.1).  That is, the following conjecture
may still be true.

\begin{conjecture}\label{conjecture2}
Let $I=(-1, 1)$ and let $n$ be a positive integer. If
$t_i=\frac{it}{n}$ for $1\leq i \leq n$, then the function
\begin{equation}
F(x)=P_x\left\{\,X_{t_1}^\alpha\in I, \dots, X_{t_n}^\alpha \in
I\right\}
\end{equation}
is concave on $I$.
\end{conjecture}

A natural question is whether $\varphi_1^\alpha$ is {\it
mid--concave} for any symmetric, convex domain in the plane. We
will now show that for a large enough rhombus and $\alpha = 2$
(Brownian motion),  this is not the case. Below we use
$\lambda_1(D)$ and  $\varphi_D$  to denote the first eigenvalue
for the domain $D$ and  its corresponding eigenfunction,
respectively, for the Brownian motion. We also denote the first
exit time of the Brownian motion from a domain  $D$ by  $\tau_D$.

\begin{proposition}\label{proposition3}
For $n\geq 1$, set $$ D(n) =\left\{ \,(x_1,x_2)\in\Rt:\, x_1\in
(-n,n),\, x_2 \in \left(-1 +{\frac{
|x_1|}{n}},1-{\frac{|x_1|}{n}}\right)\, \right\}.$$ There exists
an $n$ large enough such that $\varphi_{D(n)}$ is not {\it
mid--concave} on $D(n)$.
\end {proposition}
\begin{proof}
The rectangle $$ R(n) = (-\sqrt{n},\sqrt{n}) \times \left( -1 +
{\frac{1}{\sqrt{n}}},1 - {\frac{1}{\sqrt{n}}}\,\right)$$ is a
subset of $D(n)$ and therefore, $$ \lambda_1(D(n)) <
\lambda_1(R(n)) = \frac{\pi^2}{(2 - 2/\sqrt{n})^2} +
\frac{\pi^2}{(2\sqrt{n})^2} \le \frac{\pi^2}{4} \left(1 +
\frac{3}{\sqrt{n}}\right), $$
 for $n$ large enough.
Now, for any $a \in (0,1/2)$,  consider the subset of $D(n)$ define by
\[
 Q(a,n) = \left\{\,(x_1,x_2) \in \Rt: \, x_1 \in
(an,n),  \, x_2 \in \left( -1
+{\frac{|x_1|}{n}},1-{\frac{|x_1|}{n}}\right)\,\right\}.
\]
Since $$ \Delta \varphi_{D(n)} +\lambda_1(D(n))\varphi_{D(n)}=0 $$
in $Q(a,n)$ and $\lambda_1(D(n))< \lambda_1(Q(a,n))$. That is,
$\varphi_{D(n)}$ is a $q$-harmonic function with $q =
\lambda_1(D(n))$.  The Feynman--Kac formula gives that  for any $x
\in Q(a,n)$,
\begin{eqnarray}
\varphi_{D(n)}(x) &=& E_x \left[\, e^{\lambda_1(D(n))\,
\tau_{Q(a,n)}}\, \;\varphi_{D(n)}(\,B(\tau_{Q(a,n)}))\;
\right]\nonumber\\ \label{Holder} & \le &\,\varphi_{D(n)}(0) \;
 E_x \left[\,e^{\lambda_1(D(n)) \,\tau_{Q(a,n)}}
; \, B(\tau_{Q(a,n)}) \in D(n) \setminus Q(a,n) \,\right].
\end{eqnarray}
Of course, $$\varphi_{D(n)}(0) = \max\left\{\varphi_{D(n)}(x):  x
\in D\right\},$$ by symmetry.
 Let $p(a) = (1 - a/2)^2/(1 - a)^2$ and $q(a)$ be such that
$1/p(a) + 1/q(a) = 1$. Note that $p(a) > 1$ so $q(a) > 0$. By
H{\"o}lder's  inequality the expression in (\ref{Holder}) is bounded
 above by
\begin{eqnarray*}
&& \varphi_{D(n)}(0) \;\left(\,E_x \left[\,e^{\lambda_1(D(n))\,
\tau_{Q(a,n)}\, p(a)}\,\right]\,\right)^{1/p(a)}\\ && \times
\left(\,P_x\left[\,B(\,\tau_{Q(a,n)}) \in D(n) \setminus
Q(a,n)\,\right]\,\right)^{1/q(a)}.
\end{eqnarray*}
Since
$$Q(a,n) \subset (-\infty,\infty) \times (-1 + a, 1 - a)
$$
we have that
for any $x \in Q(a,n)$,
\begin{eqnarray*}
&& E_x \left[\, e^{\lambda_1(D(n)) \,\tau_{Q(a,n)} \,p(a)}
\,\right]\\ && \le E_{0} \left[\, \exp \left( (\pi^2/4)\, (1 +
3/\sqrt{n}) \, p(a) \, \tau_{(-1 + a, 1 - a)}\,\right)\, \right]
\\
&& = E_{0} \left[\,\exp \left( \,(\pi^2/4)\, (1 + 3/\sqrt{n}) \,
 p(a)\, (1 - a)^2 \, \tau_{(-1, 1)}\,\right) \,\right]
 \\
&&= E_{0} \left[\,\exp \left(\, (\pi^2/4)\, (1 + 3/\sqrt{n}) \,
  (1 - a/2)^2 \,  \tau_{(-1, 1)}\,\right) \,\right].
\end{eqnarray*}
By a simple calculation we see  that $(1 + 3/\sqrt{n}) (1 - a/2)
\le 1$ when $n \ge (6 - 3a)^2/a^2$. For such $n$, we have
\begin{eqnarray*}
&& \left( \,E_x \left[ \,e^{\,\lambda_1(D(n)) \,\tau_{Q(a,n)}\,
p(a)} \,\right]\, \right)^{1/p(a)}
\\
&&
\le
\left(\,E_{0} \left[\,\exp\left( \,(\pi^2/4) \,
  (1 - a/2) \, \tau_{(-1, 1)}\right)\, \right]\, \right)^{1/p(a)}
= C_1(a).
\end{eqnarray*}
Using the fact that $\frac{\pi^2}{4}$ is the eigenvalue for the
interval $(-1, 1)$, we have that for any $c \in (0,1)$,
$E_0[\,\exp(c\,\tau_{(-1,1)}\,  \pi^2/4)\,] < \infty$. Thus  for
any $a \in (0,1/2)$ we have $C_1(a) < \infty$.

By standard results for Brownian motion (or the trivial estimate
of the harmonic measure in the strip obtained by conformal mapping
to the disk),  for any $b \ge 0$ and $x_1 > b$ we have $$
P_{(x_1,0)}\left[ \, B(\tau_{(b,\infty) \times (-1,1)}) \in
(-\infty,b) \times (-1,1)\,\right] \le C_2 \,e^{-{\frac{\pi}{2}}
(x_1 - b)}, $$ where $C_2 > 0$ is an absolute constant.

Note that $x = (2 a n,0) \in Q(a,n)$. It follows that $$ P_{(2 a
n,0)}\left[ \,B(\tau_{Q(a,n)}) \in D(n) \setminus Q(a,n)\,\right]
\le C_2 e^{- {\frac{\pi}{2}}a n}. $$ Now choose $a = 1/8$. For
such $a$ we have $(2 a n,0) = (n/4,0)$. For $n \ge (6 - 3a)^2/a^2$
we have
\begin{equation}
\label{notconcave} \varphi_{D(n)}(n/4,0) \,\le\,
\varphi_{D(n)}(0,0) \,C_1(1/8) \, \left[ \,C_2
e^{-{\frac{\pi}{16}} n}\,\right]^{\frac{1}{q(1/8)}}.
\end{equation}

If $\varphi_{D(n)}$ were  {\it mid--concave}, we would have $$
\varphi_{D(n)}(n/4,0) \ge {\frac{1}{2}}\,[\,\varphi_{D(n)}(0,0) +
\varphi_{D(n)}(n/2,0)\,] \ge {\frac{1}{2}}\varphi_{D(n)}(0,0). $$
However, by (\ref{notconcave}) for large enough $n$  we have that
$\varphi_{D(n)}(n/4,0)$ is smaller than $\varphi_{D(n)}(0,0)/2 .$
Thus $\varphi_{D(n)}$ is not {\it mid--concave.} Indeed, the same
argument shows that  for any $c \in (0,1)$ there exists an $n$
large enough such that $\varphi_{D(n)}$ is not concave on the
interval with endpoints $(-cn,0)$, $(cn,0)$.  \end{proof}

\end{document}